\documentclass[a4paper,12pt]{article}

\usepackage{amsmath,amssymb,amscd,cite}

\makeatletter
 
  \@addtoreset{equation}{section}
 \makeatother

\newtheorem{definition}{{\bf Definition}}[section]
\newtheorem{theorem}[definition]{{\bf Theorem}}
\newtheorem{lemma}[definition]{{\bf Lemma}}
\newtheorem{proposition}[definition]{{\bf Proposition}}

\begin{document}

\begin{center}

\textbf{\Large Non-commutative hypergroup of order five}

\bigskip

{\bf Yasumichi Matsuzawa, Hiromichi Ohno, Akito Suzuki, \\
Tatsuya Tsurii and Satoe Yamanaka}

\bigskip

\end{center}

\begin{abstract}

We prove that all hypergroups of order four are commutative
and that there exists a non-comutative hypergroup of order five.
These facts imply that 
the minimum order of non-commutative hypergroups is five
even though the minimum order of non-commutative groups is six.

\bigskip

\medskip



\end{abstract}


\section{Introduction}

\medskip
The structure of hypergroups has been studied by many authors
(see
\cite{BH,HJKK,HK1,HK2,HKKK,IK1,IK2,IKS,K,KKY,KM,KMTY,KST,KSTY,
KSY,KY,W1,W2} and references therein).
N. J. Wildberger \cite{W2} determined
the structure of hypergroups of order three;
however the structure of hypergroups of order greater than three
has not been determined. 
R. Ichihara and S. Kawakami \cite{IK1}
revealed the structure of hypergroups of order four
which has a subhypergroup.
However, even the structure of a hypergroup of order four
 has not been determined in general. 

In this paper, we focus on the commutativity of finite hypergroups. 
Let $\mathcal{K} = \{c_0, c_1, \dots, c_n\}$ be a finite hypergroup of order $n$.
We first prove that $\mathcal{K}$ is commutative
if the $*$-structure of $\mathcal{K}$ is trivial,
{\it i.e.}, $c_i^* = c_i$ ($i=0, 1, \dots, n$).
From this fact, all hypergroups of order less than three are commutative.
As was shown in \cite{W2}, 
all hypergroups of order three are commutative.
We next  prove that all hypergroups of order four are commutative.
Therefore, all hypergroups of order less than five are commutative.
Let us consider the case where $\mathcal{K} = \{c_0, c_1, c_2, c_3, c_4\}$
is a hypergroup of order five.
Then, the following three $*$-structures are possible:

\medskip
(1)~$c_0^* = c_0$, $c_1^* = c_1$, $c_2^* = c_2$, $c_3^* = c_3$, $c_4^* = c_4$ (trivial case);

\medskip

(2)~$c_0^* = c_0$, $c_1^* = c_4$, $c_2^* = c_3$, $c_3^* = c_2$, $c_4^* = c_1$;

\medskip

(3)~$c_0^* = c_0$, $c_1^* = c_1$, $c_2^* = c_2$, $c_3^* = c_4$, $c_4^* = c_3$.

\medskip
\noindent
In case (1),  the $*$-structure is trivial and hence $\mathcal{K}$ is commutative.
We prove that in case (2), $\mathcal{K}$ is commutative,
but case (3) allows $\mathcal{K}$ to be non-commutative. 
We emphasize that
the minimum order of non-commutative hypergroups is five,
whereas that of non-commutative groups is six.

This paper is organized as follows.
In Section 2, we review some of the standard facts on finite hypergroups.
In Section 3, we prove that 
all hypergroups of order less than five
and hypergroups with a trivial $*$-structure are commutative.
Section 4 is devoted to the study of hypergroups of order five.
In the appendix, we present the detailed calculations of associativities.

\bigskip


\section{Preliminaries}

\medskip
For a finite set $K = \{c_0, c_1, \dots, c_n\}$, we denote
the set of all complex-valued measures on $K$ and the set of all non-negative 
probability measures on $K$ by $M^b(K)$ and $M^1(K)$, respectively:
\begin{align*}
	\begin{split}
M^b(K) &:= \left\{\sum_{j = 0}^n a_j\delta_{c_j}~:~a_j \in \mathbb{C}~~
(j = 0, 1, 2, \dots, n)\right\}, \\
M^1(K) &:= \left\{\sum_{j = 0}^n a_j\delta_{c_j}~:~a_j \geq 0 ~~
(j = 0,1, 2, \dots, n), ~\sum_{j=0}^n a_j = 1\right\},
    \end{split}
\end{align*}
where the symbol $\delta_c$ stands for the Dirac measure at $c \in K$. 
The support of $\mu = a_0\delta_{c_0} + a_1\delta_{c_1} + \cdots + a_n\delta_{c_n} \in M^b(K)$ is
\[
\text{supp}(\mu) := \{c_j \in K~:~a_j \not = 0~(j = 0,1, 2, \dots, n)\}. 
\]

\medskip

\noindent{\bf Axiom.}  
A {\it finite  hypergroup} $\mathcal{K}= (K, M^b(K), \circ, *)$ comprises a 
finite set $K =\{c_0, c_1, \dots, c_n\}$ 
together with an associative product (called convolution) $\circ$
and an involution $\ast$ on $M^b(K)$, satisfying the following conditions. 

\medskip
\begin{itemize}
\item[(a)] The space ($M^b(K), \circ, \ast$)  is an associative $\ast$-algebra 
with unit $\delta_{c_0}$;
\item[(b)] For $c_i, c_j \in K$, the convolution $\delta_{c_i} \circ \delta_{c_j}$ 
belongs to $M^1(K)$;
\item[(c)] There exists an involutive bijection $c_i \mapsto c_i^\ast$ on $K$ such that 
$\delta_{c_i^\ast} = \delta_{c_i}^\ast$.  
Moreover, $c_j = c_i^\ast$ if and only if  
$c_0 \in \text{supp}(\delta_{c_i} \circ \delta_{c_j})$ for all $c_i, c_j \in K$.  
\end{itemize}

\noindent A finite hypergroup $\mathcal{K}$ is called {\it commutative} if 
the convolution $\circ$ on $M^b(K)$ is commutative. 
Sometimes, we express the hypergroup $\mathcal{K}= (K, M^b(K), \circ, *)$ in  
$\mathcal{K} = \{c_0,c_1, \dots, c_n\}$. 
\medskip

From axiom (c), the structure equation of $\delta_{c_i} \circ \delta_{c_i^*}$ is written as 
\[
\delta_{c_i} \circ \delta_{c_i^*} = 
a_0 \delta_{c_0} + a_1 \delta_{c_1} + \cdots + a_n\delta_{c_n}
\quad(a_0 > 0).
\]
Then, the weight of an element $c_i \in K$ is defined as $w(c_i) := \frac{1}{a_0}$, 
and the total weight of $\mathcal{K}$ is given by 
$w(\mathcal{K}):= \sum_{i=0}^n w(c_i)$. 
We use $\delta_{c_i} \circ \delta_{c_j}|_{c_0}$ 
to denote the coefficient of 
$\delta_{c_0}$ of the convolution $\delta_{c_i} \circ \delta_{c_j}$. 


\medskip
\begin{lemma} Let $\mathcal{K}$ be a finite hypergroup. For $c_i \in K$, 
\begin{equation}\label{eq1}
\delta_{c_i} \circ \delta_{c_i^*}|_{c_0} =  
\delta_{c_i^*} \circ \delta_{c_i}|_{c_0} .
\end{equation}
\end{lemma}

\noindent
{\bf Proof.}
The lemma follows from the fact that 
\begin{equation}
\label{eq01}
 w(c_i) = w(c_i^*), \quad c_i \in K
 \end{equation}
(see \cite[Lemma 1.3]{SW1}). 
\hspace{\fill}$\Box$

\bigskip

We write $\delta_{c_i} \circ \delta_{c_j}$ as
\[
\delta_{c_i} \circ \delta_{c_j} = \sum_{k=0}^n m_{ij}^k \delta_{c_k}.
\]

\begin{lemma}\label{lemma2.2}
For all $1\le i,j,k \le n$,
\begin{eqnarray}
m_{ij}^k = m_{\bar{j}\bar{i}}^{\bar{k}}, \label{eq2} \\
m_{ij}^{\bar{i}} = m_{ji}^{\bar{i}}, \label{eq3} \\
m_{ij}^i = m_{j\bar{i}}^{\bar{i}}, \label{eq4} 
\end{eqnarray}
where $\bar{i}$ satisfies $c_i^* = c_{\bar{i}}$.
\end{lemma}

\noindent
{\bf Proof.}
By the $*$-structure,
$ (\delta_{c_i} \circ \delta_{c_j} )^* = \delta_{c_j}^* \circ \delta_{c_i}^*$.
Then, \eqref{eq2} follows from
\begin{eqnarray*}
(\delta_{c_i} \circ \delta_{c_j} )^* =  \left(\sum_{k=0}^n m_{ij}^k \delta_{c_k} \right)^*
 =  \sum_{k=0}^n m_{ij}^k \delta_{c_k}^*, \\
\delta_{c_j}^* \circ \delta_{c_i}^* =
\delta_{c_j^*} \circ \delta_{c_i^*} =\delta_{c_{\bar{j}}} \circ \delta_{c_{\bar{i}}} =
  \sum_{k=0}^n m_{\bar{j}\bar{i}}^k \delta_{c_k}.
\end{eqnarray*}
 
The coefficients $\delta_{c_i} \circ (\delta_{c_j}\circ \delta_{c_i}) |_{c_0}$ and
$(\delta_{c_i} \circ \delta_{c_j}) \circ \delta_{c_i} |_{c_0}$ are calculated as
\begin{eqnarray*}
\delta_{c_i} \circ (\delta_{c_j}\circ \delta_{c_i}) |_{c_0} 
=\delta_{c_i} \circ ( \sum_{k=0}^n m_{ji}^k \delta_{c_k}) |_{c_0} = w(c_i) m_{ji}^{\bar{i}}, \\
(\delta_{c_i} \circ \delta_{c_j}) \circ \delta_{c_i} |_{c_0}
=(\sum_{k=0}^n m_{ij}^k \delta_{c_k}) \delta_{c_i} |_{c_0} = w(c_i^*) m_{ij}^{\bar{i}}.
\end{eqnarray*}
By \eqref{eq01} and the associativity law,
we obtain \eqref{eq3}. 

Similarly, the equations  
\begin{eqnarray*}
\delta_{c_i} \circ (\delta_{c_j}\circ \delta_{c_i}^*) |_{c_0} 
=\delta_{c_i} \circ ( \sum_{k=0}^n m_{j\bar{i}}^k \delta_{c_k}) |_{c_0} 
= w(c_i) m_{j\bar{i}}^{\bar{i}}, \\
(\delta_{c_i} \circ \delta_{c_j}) \circ \delta_{c_i}^* |_{c_0}
=(\sum_{k=0}^n m_{ij}^k \delta_{c_k}) \circ \delta_{c_i}^*|_{c_0} = w(c_i) m_{ij}^{i}
\end{eqnarray*}
yield \eqref{eq4}, which completes the proof.
\hspace{\fill}$\Box$


\section{Commutativity of hypergroups}

In this section, we discuss the commutativity of hypergroups. First, we prove the
following proposition.

\begin{proposition}\label{prop2}
Let $\mathcal{K}$ be a finite hypergroup. If the $*$-structure of $\mathcal{K}$ is trivial, 
{\it i.e.}, for all $c_i \in K$,  $c_i^* = c_i$, then $\mathcal{K}$ is commutative. 
\end{proposition}

\noindent
{\bf Proof.} Let $\delta_{c_i} \circ \delta_{c_j} = \sum_{k=0}^n m_{ij}^k \delta_{c_k}$ 
where $\sum_{k=0}^n m_{ij}^k = 1$ and $m_{ij}^k \geq 0$. 
By \eqref{eq2} and the triviality of the $*$-structure,
we have $m_{ij}^k = m_{\bar{j}\bar{i}}^{\bar{k}} = m_{ji}^k$.
Hence, for all $1\le i,j \le n$,
\[
\delta_{c_i} \circ \delta_{c_j} 
=\sum_{k=0}^n m_{ij}^k \delta_{c_k}
=\sum_{k=0}^n m_{ji}^k \delta_{c_k}
= \delta_{c_j} \circ \delta_{c_i}.
\]
When $i=0$ or $j=0$, the commutativity of $\delta_{c_i} $ and 
$\delta_{c_j}$ is trivial. 
\hspace{\fill}$\Box$

\bigskip

Proposition \ref{prop2} yields that 
all hypergroups of order less than three are commutative.
The commutativity of a hypergroup of order three is already known.

\begin{proposition}{\rm \cite{W2}}
 Let $\mathcal{K}$ be a hypergroup of order three. Then, $\mathcal{K}$ is commutative. 
\end{proposition}

All hypergroups of order four are also commutative.

\begin{theorem}
 Let $\mathcal{K}$ be a hypergroup of order four. Then, $\mathcal{K}$ is  commutative. 
\end{theorem}

\noindent
{\bf Proof.} Let $\mathcal{K} = \{c_0, c_1, c_2, c_3\}$ be a hypergroup of order four. 
All possible $*$-structures of $\mathcal{K}$ are as follows:

\medskip
(i)~$c_0^* = c_0$, $c_1^* = c_1$, $c_2^* = c_2$, $c_3^* = c_3$ (trivial case);

(ii)~$c_0^* = c_0$, $c_1^* = c_3$, $c_2^* = c_2$, $c_3^* = c_1$. 

\medskip
\noindent
Note that the case in which $c_0^* = c_0$, $c_1^* = c_2$, $c_2^* = c_1$, $c_3^* = c_3$ 
is essentially the same as case (ii). 
By Proposition \ref{prop2}, $\mathcal{K}$ is commutative in case (i); thus,
we need to consider only case (ii). 
Although the proof is similar to that of Lemma \ref{lemma2.2},
we provide the proof for completeness.

By axiom (c) of a finite hypergroup, 
the structure equation of 
$\delta_{c_1} \circ \delta_{c_3}$ is written as
$$
\delta_{c_1} \circ \delta_{c_3} = p_0 \delta_{c_0} + p_1 \delta_{c_1} + 
p_2 \delta_{c_2} + p_3 \delta_{c_3}~~~(p_0 > 0). 
$$
Combining this with $c_1^*=c_3$ yields
\begin{eqnarray*}
(\delta_{c_1} \circ \delta_{c_3})^* &=& (p_0 \delta_{c_0} + p_1 \delta_{c_1} + 
p_2 \delta_{c_2} + p_3 \delta_{c_3})^*, \\
\delta_{c_1} \circ \delta_{c_3} &=& p_0 \delta_{c_0} + p_1 \delta_{c_3} + 
p_2 \delta_{c_2} + p_3 \delta_{c_1}.
\end{eqnarray*}
Hence, we have $p_1 = p_3$ and
\begin{equation}
\label{eq03}
\delta_{c_1} \circ \delta_{c_3} = p_0 \delta_{c_0} + p_1 \delta_{c_1} + 
p_2 \delta_{c_2} + p_1 \delta_{c_3}. 
\end{equation}
By axiom (c) and \eqref{eq1},
the structure equation of $\delta_{c_3} \circ \delta_{c_1}$ is  
$
\delta_{c_3} \circ \delta_{c_1} 
= p_0 \delta_{c_0} + q_1 \delta_{c_1} + q_2 \delta_{c_2} + q_3 \delta_{c_3}$ ($p_0 > 0$).
Similar to the above mentioned instance, we know that $q_1= q_3$ and
\begin{equation}
\label{eq04}
\delta_{c_3} \circ \delta_{c_1} 
= p_0 \delta_{c_0} + q_1 \delta_{c_1} + q_2 \delta_{c_2} + q_1 \delta_{c_3}. 
\end{equation}
The structure equation of $\delta_{c_2} \circ \delta_{c_2}$ is written as
\begin{align*}
\delta_{c_2} \circ \delta_{c_2} 
= r_0 \delta_{c_0} + r_1 \delta_{c_1} + r_2 \delta_{c_2} + r_1 \delta_{c_3}~~~(r_0 > 0). 
\end{align*}
By axiom (c), we can write
$$
\delta_{c_1} \circ \delta_{c_1} 
= \alpha_1 \delta_{c_1} + \alpha_2 \delta_{c_2} + \alpha_3 \delta_{c_3}. 
$$
Hence, by the $*$-structure 
$(\delta_{c_1} \circ \delta_{c_1})^* 
	= \delta_{c_3} \circ \delta_{c_3}$, 
we obtain
\[
\delta_{c_3} \circ \delta_{c_3} 
= \alpha_3 \delta_{c_1} + \alpha_2 \delta_{c_2} + \alpha_1 \delta_{c_3}.
\]
Similarly, we can write
\begin{equation}
\label{eq09}
\begin{cases}
\delta_{c_1} \circ \delta_{c_2} 
= \beta_1 \delta_{c_1} + \beta_2 \delta_{c_2} + \beta_3 \delta_{c_3},\\
\delta_{c_2} \circ \delta_{c_3} 
= \beta_3 \delta_{c_1} + \beta_2 \delta_{c_2} + \beta_1 \delta_{c_3},\\
\delta_{c_2} \circ \delta_{c_1} 
= \gamma_1 \delta_{c_1} + \gamma_2 \delta_{c_2} + \gamma_3 \delta_{c_3},\\
\delta_{c_3} \circ \delta_{c_2} 
= \gamma_3 \delta_{c_1} + \gamma_2 \delta_{c_2} + \gamma_1 \delta_{c_3}. 
\end{cases}
\end{equation}

Next, we consider the associativity of $\mathcal{K}$. 
By axiom (a),
$$
(\delta_{c_3} \circ \delta_{c_1} ) \circ \delta_{c_1}|_{c_0} =
\delta_{c_3} \circ (\delta_{c_1} \circ \delta_{c_1})|_{c_0},
$$
and therefore, $q_1 p_0 = \alpha_1 p_0$.  
Hence, by $p_0 > 0$, 
$$
q_1=\alpha_1.
$$
By 
$
(\delta_{c_1} \circ \delta_{c_1} ) \circ \delta_{c_3}|_{c_0} = 
\delta_{c_1} \circ (\delta_{c_1} \circ \delta_{c_3})|_{c_0} 
$ and $p_0>0$, 
we have $\alpha_1 = p_1$ and 
\begin{equation}
\label{eq02}
p_1 = q_1.
\end{equation}
By axiom (b), 
$\delta_{c_1} \circ \delta_{c_3}, \delta_{c_3} \circ \delta_{c_1}
\in M^1(K)$, with the result that \eqref{eq03} and \eqref{eq04}
yield
\begin{align*}
p_0 + 2 p_1 + p_2 = 1~~{\rm and}~~p_0 + 2q_1 + q_2 = 1. 
\end{align*}
Combining this with \eqref{eq02} 
yields
\begin{equation}
\label{eq05}
p_2 = q_2.
\end{equation}
Substituting \eqref{eq02} and \eqref{eq05} into
\eqref{eq03} and \eqref{eq04},
we conclude that
$$
\delta_{c_1} \circ \delta_{c_3} = \delta_{c_3} \circ \delta_{c_1}. 
$$ 

The associativity,  
$
(\delta_{c_2} \circ \delta_{c_1} ) \circ \delta_{c_2}|_{c_0} 
=\delta_{c_2} \circ (\delta_{c_1} \circ \delta_{c_2})|_{c_0},
$
yields $\gamma_2 r_0 = \beta_2 r_0$. 
Hence, by $r_0 > 0$,
\begin{equation}
\label{eq06}
\gamma_2 = \beta_2.
\end{equation}
Similarly, by 
$(\delta_{c_1} \circ \delta_{c_1} ) \circ \delta_{c_2}|_{c_0} =
\delta_{c_1} \circ (\delta_{c_1} \circ \delta_{c_2})|_{c_0}$
and
$(\delta_{c_2} \circ \delta_{c_1} ) \circ \delta_{c_1}|_{c_0} = 
\delta_{c_2} \circ (\delta_{c_1} \circ \delta_{c_1})|_{c_0}$,
we obtain
\begin{equation*}
\alpha_2 r_0 = \beta_3 p_0
%
%
\quad \mbox{and} \quad
\gamma_3 p_0 = \alpha_2 r_0.
\end{equation*}
Combining this with $p_0 > 0$
yields
\begin{equation}
\label{eq07}
\beta_3 = \gamma_3.
\end{equation}
Because
$\delta_{c_1} \circ \delta_{c_2}, \delta_{c_2} \circ \delta_{c_1}
\in M^1(K)$,
\begin{equation*}
\beta_1 + \beta_2 + \beta_3 = 1~~{\rm and}~~\gamma_1 + \gamma_2 + \gamma_3 = 1. 
\end{equation*}
Hence, by \eqref{eq06} 
and  \eqref{eq07}, 
we obtain
\begin{equation}
\label{eq08}
\beta_1 = \gamma_1.
\end{equation}
Substituting \eqref{eq06}, \eqref{eq07}, and \eqref{eq08} 
into \eqref{eq09}
yields
$$
\delta_{c_1} \circ \delta_{c_2} = \delta_{c_2} \circ \delta_{c_1}~~{\rm and}~~
\delta_{c_2} \circ \delta_{c_3} = \delta_{c_3} \circ \delta_{c_2}.
$$
This completes the proof.
\hspace{\fill}$\Box$


\section{Non-commutative hypergroup of order five}

\medskip
In this section, we discuss the commutativity and non-commutativity of hypergroups 
of order five. Consequently, we can find a non-commutative hypergroup of 
order five. 
Let $\mathcal{K} = \{c_0, c_1, c_2, c_3, c_4\}$ be a hypergroup of order five. 
All possible $*$-structures of $\mathcal{K}$ are as follows:
\begin{itemize}
\item[(1)] $c_0^* = c_0$, $c_1^* = c_1$, $c_2^* = c_2$, $c_3^* = c_3$, $c_4^* = c_4$ (trivial case);
\item[(2)] $c_0^* = c_0$, $c_1^* = c_4$, $c_2^* = c_3$, $c_3^* = c_2$, $c_4^* = c_1$;
\item[(3)] $c_0^* = c_0$, $c_1^* = c_1$, $c_2^* = c_2$, $c_3^* = c_4$, $c_4^* = c_3$.
\end{itemize}

\noindent
By Proposition \ref{prop2},
$\mathcal{K}$ is commutative  in case (1).
Hence, we need to consider only cases (2) and (3).

\begin{proposition}
\label{prop4.1}
 Let $\mathcal{K} = \{c_0, c_1, c_2, c_3, c_4\}$ be a hypergroup of 
order five. 
If  $\mathcal{K}$ satisfies (2), 
then $\mathcal{K}$ is commutative. 
\end{proposition}

\medskip
\noindent
{\bf Proof.} 
If we write 
$\delta_{c_i} \circ \delta_{c_j} 
	= \sum_{k=0}^n a_k \delta_{c_k}$,
then
$\delta_{c_j}^* \circ \delta_{c_i}^* 
	= \sum_{k=0}^n a_k \delta_{c_k}^*$,
because $\delta_{c_j}^* \circ \delta_{c_i}^* 
	=(\delta_{c_i} \circ \delta_{c_j} )^*$.
By axiom (c),
$\delta_{c_i} \circ \delta_{c_j}|_{c_0} > 0$ only when $c_j = c_i^*$.
Moreover, by \eqref{eq1}, 
$\delta_{c_i} \circ \delta_{c_i}^*|_{c_0} = \delta_{c_i}^* \circ \delta_{c_i}|_{c_0}$.
Combining these,
we can write the structure equations of $\mathcal{K}$ as  
\begin{align*}
\delta_{c_1} \circ \delta_{c_1} &= 
a_1 \delta_{c_1} + a_2 \delta_{c_2} + a_3 \delta_{c_3} + a_4 \delta_{c_4},\\
\delta_{c_4} \circ \delta_{c_4} &= 
a_4 \delta_{c_1} + a_3 \delta_{c_2} + a_2 \delta_{c_3} + a_1 \delta_{c_4}, 
\end{align*}
\begin{align*}
\delta_{c_2} \circ \delta_{c_2} &= 
b_1 \delta_{c_1} + b_2 \delta_{c_2} + b_3 \delta_{c_3} + b_4 \delta_{c_4},\\
\delta_{c_3} \circ \delta_{c_3} &= 
b_4 \delta_{c_1} + b_3 \delta_{c_2} + b_2 \delta_{c_3} + b_1 \delta_{c_4},
\end{align*}
\begin{align*}
\delta_{c_1} \circ \delta_{c_2} &=  
d_1 \delta_{c_1} + d_2 \delta_{c_2} + d_3 \delta_{c_3} + d_4 \delta_{c_4},\\
\delta_{c_3} \circ \delta_{c_4} &=
d_4 \delta_{c_1} + d_3 \delta_{c_2} + d_2 \delta_{c_3} + d_1 \delta_{c_4},\\
\delta_{c_2} \circ \delta_{c_1} &= 
e_1 \delta_{c_1} + e_2 \delta_{c_2} + e_3 \delta_{c_3} + e_4 \delta_{c_4},\\
\delta_{c_4} \circ \delta_{c_3} &=
e_4 \delta_{c_1} + e_3 \delta_{c_2} + e_2 \delta_{c_3} + e_1 \delta_{c_4},
\end{align*}
\begin{align*}
\delta_{c_1} \circ \delta_{c_3} &=
f_1 \delta_{c_1} + f_2 \delta_{c_2} + f_3 \delta_{c_3} + f_4 \delta_{c_4},\\
\delta_{c_2} \circ \delta_{c_4} &=
f_4 \delta_{c_1} + f_3 \delta_{c_2} + f_2 \delta_{c_3} + f_1 \delta_{c_4},\\
\delta_{c_3} \circ \delta_{c_1} &= 
g_1 \delta_{c_1} + g_2 \delta_{c_2} + g_3 \delta_{c_3} + g_4 \delta_{c_4},\\
\delta_{c_4} \circ \delta_{c_2} &= 
g_4 \delta_{c_1} + g_3 \delta_{c_2} + g_2 \delta_{c_3} + g_1 \delta_{c_4},
\end{align*}
\begin{align*}
\delta_{c_1} \circ \delta_{c_4} &=
h_0 \delta_{c_0} + h_1 \delta_{c_1} + h_2 \delta_{c_2} + h_2 \delta_{c_3} + h_1 \delta_{c_4},\\
\delta_{c_4} \circ \delta_{c_1} &= 
h_0 \delta_{c_0} + i_1 \delta_{c_1} + i_2 \delta_{c_2} + i_2 \delta_{c_3} + i_1 \delta_{c_4},
\end{align*}
\begin{align*}
\delta_{c_2} \circ \delta_{c_3} &=
j_0 \delta_{c_0} + j_1 \delta_{c_1} + j_2 \delta_{c_2} + j_2 \delta_{c_3} + j_1 \delta_{c_4},\\
\delta_{c_3} \circ \delta_{c_2} &=
j_0 \delta_{c_0} + k_1 \delta_{c_1} + k_2 \delta_{c_2} + k_2 \delta_{c_3} + k_1 \delta_{c_4},
\end{align*}
where $h_0 > 0$ and $j_0 > 0$. 
Note that each block in this list contains 
$\delta_{c_i} \circ \delta_{c_j}$, $\delta_{c_j}^* \circ \delta_{c_i^*}$, $\delta_{c_j} \circ \delta_{c_i}$,
and $\delta_{c_i}^* \circ \delta_{c_j}^*$.
To prove commutativity, 
we only need to show that 
$ \delta_{c_i} \circ \delta_{c_j} = \delta_{c_j} \circ \delta_{c_i}$ 
for all $1\le i,j \le 4$ with $i \neq j$.

By \eqref{eq3}, 
we get 
\[
i_1 = m_{41}^1 = m_{14}^1 = h_1 \quad (i=4, j= 1).
\]
By axiom (b), 
$h_0 + 2h_1+2h_2 = h_0 + 2i_1 +2i_2 =1$. 
Therefore, $h_2 = i_2$ and 
$\delta_{c_1} \circ \delta_{c_4} = \delta_{c_4} \circ \delta_{c_1}$.

Similarly, we observe from \eqref{eq3} that
\[
j_2 = m_{23}^3 = m_{32}^3 = k_2 \quad (i=2, j=3).
\]
By axiom (b), $j_1 = k_1$,
and hence,
$\delta_{c_2} \circ \delta_{c_3} = \delta_{c_3} \circ \delta_{c_2}$.

Next, we consider the coefficients of $\delta_{c_1} \circ \delta_{c_2}$
and $\delta_{c_2} \circ \delta_{c_1} $.
Again by \eqref{eq3}, we have
\begin{align*}
& e_3 = m_{21}^3 = m_{12}^3 = d_3 \quad (i=2, j=1) \\
& d_4 = m_{12}^4 = m_{21}^4 = e_4 \quad (i=1, j=2).
\end{align*}
From the associativity law,
\[
(\delta_{c_1} \circ \delta_{c_2}) \circ \delta_{c_3}|_{c_0} = 
\delta_{c_1} \circ (\delta_{c_2} \circ \delta_{c_3})|_{c_0},~~~
(\delta_{c_3} \circ \delta_{c_2}) \circ \delta_{c_1}|_{c_0} = 
\delta_{c_3} \circ (\delta_{c_2} \circ \delta_{c_1})|_{c_0},
\]
and thus, we obtain
\[
d_2 j_0 = h_0 j_1 \quad \text{and} \quad k_1 h_0 = e_2 j_0.
\]
Since $j_1 = k_1$ and $j_0>0$, we get $d_2 = e_2$.
By axiom (b), $d_1 = e_1$; hence,
$\delta_{c_1} \circ \delta_{c_2} = \delta_{c_2} \circ \delta_{c_1}$ and
$\delta_{c_3} \circ \delta_{c_4} = \delta_{c_4} \circ \delta_{c_3}$.

Finally, we consider the coefficients of $\delta_{c_1} \circ \delta_{c_3}$
and $\delta_{c_3} \circ \delta_{c_1} $.
By \eqref{eq3}, 
\begin{align*}
& f_2 = m_{24}^3 = m_{42}^3 = g_2 \quad (i=2, j=4)\\
& f_4 = m_{13}^4 = m_{31}^4 = g_4 \quad (i=1, j=3).
\end{align*}
By \eqref{eq4}, 
\begin{align*}
& f_1 = m_{13}^1 = m_{34}^4 =d_1 \quad (i=1, j=3), \\
& e_1 = m_{43}^4 = m_{31}^1 = g_1 \quad (i=4, j=3).
\end{align*}
Since $e_1 =d_1$, $f_1 = g_1$.
Again by axiom (b), we obtain $f_3 = g_3$.
Therefore,
$\delta_{c_1} \circ \delta_{c_3} = \delta_{c_3} \circ \delta_{c_1}$ and
$\delta_{c_2} \circ \delta_{c_4} = \delta_{c_4} \circ \delta_{c_2}$.
\hspace{\fill}$\Box$

\bigskip


In what follows, we consider a hypergroup with the $*$-structure (3). 

\begin{theorem}
\label{thm}
There exists a non-commutative hypergroup of order five 
with the $*$-structure (3).
\end{theorem}

\medskip
\noindent
In the remainder of the paper,
we prove Theorem \ref{thm}.
As in the proof of Proposition \ref{prop4.1}, 
the structure equations of $\mathcal{K}$ are 
\begin{align*}
&\delta_{c_1} \circ \delta_{c_1} 
= p_0 \delta_{c_0} + p_1\delta_{c_1} + p_2 \delta_{c_2} + p_3\delta_{c_3} + p_3\delta_{c_4},
\end{align*}
\begin{align*}
&\delta_{c_2} \circ \delta_{c_2} 
= q_0 \delta_{c_0} + q_1\delta_{c_1} + q_2 \delta_{c_2} + q_3\delta_{c_3} + q_3\delta_{c_4},
\end{align*}
\begin{align*}
&\delta_{c_3} \circ \delta_{c_4} 
= r_0 \delta_{c_0} + r_1\delta_{c_1} + r_2 \delta_{c_2} + r_3\delta_{c_3} + r_3\delta_{c_4},\\
&\delta_{c_4} \circ \delta_{c_3} 
= r_0 \delta_{c_0} + s_1\delta_{c_1} + s_2 \delta_{c_2} + s_3\delta_{c_3} + s_3\delta_{c_4},
\end{align*}
\begin{align*}
&\delta_{c_1} \circ \delta_{c_2} 
= \alpha_1\delta_{c_1} + \alpha_2 \delta_{c_2} + \alpha_3\delta_{c_3} + \alpha_4\delta_{c_4},\\
&\delta_{c_2} \circ \delta_{c_1} 
= \alpha_1\delta_{c_1} + \alpha_2 \delta_{c_2} + \alpha_4\delta_{c_3} + \alpha_3\delta_{c_4},
\end{align*}
\begin{align*}
&\delta_{c_1} \circ \delta_{c_3} 
= \beta_1\delta_{c_1} + \beta_2 \delta_{c_2} + \beta_3\delta_{c_3} + \beta_4\delta_{c_4},\\
&\delta_{c_4} \circ \delta_{c_1} 
= \beta_1\delta_{c_1} + \beta_2 \delta_{c_2} + \beta_4\delta_{c_3} + \beta_3\delta_{c_4},\\
&\delta_{c_3} \circ \delta_{c_1} 
= \gamma_1\delta_{c_1} + \gamma_2 \delta_{c_2} + \gamma_4\delta_{c_3} + \gamma_3\delta_{c_4},\\
&\delta_{c_1} \circ \delta_{c_4} 
= \gamma_1\delta_{c_1} + \gamma_2 \delta_{c_2} + \gamma_3\delta_{c_3} + \gamma_4\delta_{c_4},
\end{align*}
\begin{align*}
&\delta_{c_2} \circ \delta_{c_3} 
= \epsilon_1\delta_{c_1} + \epsilon_2 \delta_{c_2} + 
\epsilon_3\delta_{c_3} + \epsilon_4\delta_{c_4},\\
&\delta_{c_4} \circ \delta_{c_2} 
= \epsilon_1\delta_{c_1} + \epsilon_2 \delta_{c_2} + 
\epsilon_4\delta_{c_3} + \epsilon_3\delta_{c_4},\\
&\delta_{c_3} \circ \delta_{c_2} 
= \zeta_1\delta_{c_1} + \zeta_2 \delta_{c_2} + \zeta_4\delta_{c_3} + \zeta_3\delta_{c_4},\\
&\delta_{c_2} \circ \delta_{c_4} 
= \zeta_1\delta_{c_1} + \zeta_2 \delta_{c_2} + \zeta_3\delta_{c_3} + \zeta_4\delta_{c_4},
\end{align*}
\begin{align*}
&\delta_{c_3} \circ \delta_{c_3} 
= \eta_1\delta_{c_1} + \eta_2 \delta_{c_2} + \eta_3\delta_{c_3} + \eta_4\delta_{c_4},\\
&\delta_{c_4} \circ \delta_{c_4} 
= \eta_1\delta_{c_1} + \eta_2 \delta_{c_2} + \eta_4\delta_{c_3} + \eta_3\delta_{c_4},
\end{align*}
where $p_0 > 0$, $q_0 > 0$, and $r_0 > 0$. 
By \eqref{eq3}, 
 we get 
\begin{align}
& \beta_1= m_{13}^1 = m_{31}^1 = \gamma_1 \quad (i=1, j=3), 
\label{eq012} \\
& \gamma_3= m_{31}^4 = m_{13}^4 = \beta_4 \quad (i=3, j=1), 
\label{eq011} \\
& \epsilon_2= m_{23}^2 = m_{32}^2 = \zeta_2 \quad (i=2, j=3), 
\label{eq014} \\
& \zeta_3= m_{32}^4 = m_{23}^4 = \epsilon_4 \quad (i=3, j=2), 
\label{eq015} \\
& r_3= m_{34}^4 = m_{43}^4 = s_3 \quad (i=3, j=4).
\label{eq010}
\end{align}
Moreover, by \eqref{eq4}, 
we have
\begin{equation}
\label{eq016}
s_3= m_{43}^4 = m_{33}^3 = \eta_3 \quad (i=4, j=3).
\end{equation}
By axiom (b) and \eqref{eq010},
$r_1 +r_2  = s_1 +s_2 $; hence
 \begin{equation}\label{rs}
 r_1-s_1 = s_2-r_2.
\end{equation}
From the equations
\begin{align*}
&(\delta_{c_1} \circ \delta_{c_3}) \circ \delta_{c_4}|_{c_0} = 
\delta_{c_1} \circ (\delta_{c_3} \circ \delta_{c_4})|_{c_0},~~~
(\delta_{c_1} \circ \delta_{c_4}) \circ \delta_{c_3}|_{c_0} =
\delta_{c_1} \circ (\delta_{c_4} \circ \delta_{c_3})|_{c_0},
\end{align*}
we get
\begin{equation}\label{betagamma}
\beta_3 r_0 = r_1 p_0 \quad \text{and} \quad \gamma_4 r_0 = s_1 p_0.
\end{equation}
Similarly, 
\begin{align*}
(\delta_{c_2} \circ \delta_{c_3}) \circ \delta_{c_4}|_{c_0} = 
\delta_{c_2} \circ (\delta_{c_3} \circ \delta_{c_4})|_{c_0},~~~ 
(\delta_{c_2} \circ \delta_{c_4}) \circ \delta_{c_3}|_{c_0} =
\delta_{c_2} \circ (\delta_{c_4} \circ \delta_{c_3})|_{c_0},
\end{align*}
yield 
\begin{equation}\label{epsilonzeta}
\epsilon_3 r_0 = r_2 q_0 \quad \text{and} \quad \zeta_4 r_0 = s_2 q_0.
\end{equation}

\begin{lemma} 
\label{lem4.3}
If $r_1 = s_1$, then
$\mathcal{K}$ is commutative.
\end{lemma}
{\bf Proof.}
By assumption, \eqref{rs} implies that  $r_2 = s_2$. 
By axiom (b), $r_3=s_3$, and hence,
$
\delta_{c_3} \circ \delta_{c_4} = \delta_{c_4} \circ \delta_{c_3}. 
$
Similar to the above mentioned instance, 
we observe from \eqref{betagamma}  that $\beta_3 = \gamma_4$.
Combining these with axiom (b), 
we have  $\beta_2 = \gamma_2$. 
Hence,
$\delta_{c_1} \circ \delta_{c_3} = \delta_{c_3} \circ \delta_{c_1}$
and
$\delta_{c_1} \circ \delta_{c_4} = \delta_{c_4} \circ \delta_{c_1}$.
By \eqref{epsilonzeta}, 
$\epsilon_3 = \zeta_4$, and thus, $\epsilon_1 = \zeta_1$ by axiom (b).
Hence, $\delta_{c_2} \circ \delta_{c_3} = \delta_{c_3} \circ \delta_{c_2}$ and
$\delta_{c_2} \circ \delta_{c_4} = \delta_{c_4} \circ \delta_{c_2}$.
The equations
\begin{align*}
&(\delta_{c_2} \circ \delta_{c_1}) \circ \delta_{c_3}|_{c_0} = 
\delta_{c_2} \circ (\delta_{c_1} \circ \delta_{c_3})|_{c_0}~~{\rm and}~~ 
(\delta_{c_3} \circ \delta_{c_1}) \circ \delta_{c_2}|_{c_0} = 
\delta_{c_3} \circ (\delta_{c_1} \circ \delta_{c_2})|_{c_0}
\end{align*}
imply that $\alpha_3 = \alpha_4$. Hence,
$\delta_{c_1} \circ \delta_{c_2} = \delta_{c_2} \circ \delta_{c_1}.$
This completes the proof.
\hspace{\fill}$\Box$

\bigskip

To obtain a non-commutative hypergroup,
we assume that $r_1 \not = s_1$.
Then, by the structure equations,
we observe that $\delta_{c_3} \circ \delta_{c_4} \not
= \delta_{c_4} \circ \delta_{c_3}$.
Hence, $\mathcal{K}$ is non-commutative. 
\begin{lemma}
If $r_1 \not= s_1$, then
\begin{equation} 
\label{eq013}
\beta_3 \neq \gamma_4
\quad \mbox{and} \quad
\epsilon_3 \neq \zeta_4
\end{equation}
and 
\begin{align}
\alpha_1 = \beta_1 = \gamma_1,~~~\eta_4 = \beta_4 = \gamma_3 = \epsilon_4 = \zeta_3,
~~~\alpha_2 = \epsilon_2 = \zeta_2,~~{\rm and}~~\eta_3 = r_3 = s_3. 
\end{align}
\end{lemma}
{\bf Proof.}
\eqref{betagamma} and \eqref{epsilonzeta} yield \eqref{eq013}.
Now, we consider the equation
\[
(\delta_{c_3} \circ \delta_{c_4}) \circ \delta_{c_3}|_{c_4} = 
\delta_{c_3} \circ (\delta_{c_4} \circ \delta_{c_3})|_{c_4}.
\]
Combining
\begin{eqnarray*}
(\delta_{c_3} \circ \delta_{c_4}) \circ \delta_{c_3} |_{c_4}
&=& (r_0 \delta_{c_0} + r_1\delta_{c_1} + r_2 \delta_{c_2} + r_3\delta_{c_3} + r_3\delta_{c_4} )  
\circ \delta_{c_3} |_{c_4}\\
&=& r_1 \beta_4 + r_2  \epsilon_4 + r_3 (\eta_4 + s_3)
\end{eqnarray*}
and
\begin{eqnarray*}
\delta_{c_3} \circ (\delta_{c_4} \circ \delta_{c_3}) |_{c_4}
&=& \delta_{c_3}\circ
(r_0 \delta_{c_0} + s_1\delta_{c_1} + s_2 \delta_{c_2} + s_3\delta_{c_3} + s_3\delta_{c_4} )  
  |_{c_4} \\
&=& s_1 \gamma_3 + s_2  \zeta_3 + s_3 (\eta_4 + r_3)
\end{eqnarray*}
with \eqref{eq011}, \eqref{eq015}, and \eqref{eq010},
we conclude that
$\beta_4(r_1 - s_1) = \epsilon_4 (s_2 - r_2)$.
Hence, by \eqref{rs}, we get  
$$
\beta_4 = \epsilon_4.
$$ 
Similarly, calculating the coefficients
$
(\delta_{c_1} \circ \delta_{c_3}) \circ \delta_{c_1}|_{c_1}= 
\delta_{c_1} \circ (\delta_{c_3} \circ \delta_{c_1})|_{c_1}
$,
$
(\delta_{c_2} \circ \delta_{c_4}) \circ \delta_{c_2}|_{c_2} = 
\delta_{c_2} \circ (\delta_{c_4} \circ \delta_{c_2})|_{c_2}
$, 
and
$
(\delta_{c_3} \circ \delta_{c_1}) \circ \delta_{c_3}|_{c_4}= 
\delta_{c_3} \circ (\delta_{c_1} \circ \delta_{c_3})|_{c_4}
$
yields
$
\alpha_1(\gamma_2 - \beta_2) = \beta_1(\beta_3 - \gamma_4)
$,
$\alpha_2(\epsilon_1 - \zeta_1) = \epsilon_2(\zeta_4 - \epsilon_3)$,
and 
$\beta_4(\gamma_2 - \beta_2) = \eta_4(\beta_3 - \gamma_4)$. 
Combining \eqref{eq012}
--\eqref{eq015} with axiom (b) 
yields
$\gamma_2 - \beta_2 = \beta_3 - \gamma_4$
and $\epsilon_1 - \zeta_1 = \zeta_4 - \epsilon_3$.
Hence, by \eqref{eq013}, we have 
$$
\alpha_1 = \beta_1,
\quad 
\alpha_2 = \epsilon_2,
\quad \mbox{and} \quad 
\beta_4 = \eta_4. 
$$
Summarizing the above results and \eqref{eq012}--\eqref{eq016}, 
we have the lemma.
\hspace{\fill}$\Box$
%
%

\begin{lemma}
\label{lem4.5}
Let $r_1, s_1$ satisfy $r_1 \not=s_1$
and assume in addition that
\begin{equation}
\label{eq017}
p_0 = q_0 = r_0~~{\rm and}~~r_1 = s_2.
\end{equation}
Then, $r_2 = s_1$ and the following 14 equations hold:
\begin{align*}
\beta_3 = m_{13}^3 = m_{34}^{1} = r_1, \quad 
& \zeta_4 = m_{24 }^{4 }= m_{43 }^{2 } = s_2, \\ 
 \gamma_4 = m_{14 }^{4 }= m_{43}^{1 } = s_1, \quad 
& \epsilon_3 = m_{23 }^{3 }= m_{34 }^{2 } =r_2, \\ 
 p_2 = m_{11 }^{2 }= m_{12 }^{1 } = \alpha_1, \quad 
& p_3 = m_{11 }^{4 }= m_{13 }^{1 } = \beta_1, \\ 
  \alpha_2= m_{12 }^{2 }= m_{22 }^{1 } =q_1, \quad 
& q_3 = m_{22 }^{4 }= m_{23 }^{2 } = \epsilon_2, \\ 
  \beta_4 = m_{13 }^{4 }= m_{33 }^{1 } = \eta_1, \quad 
&\epsilon_4  = m_{23 }^{4}= m_{33 }^{2 } =\eta_2, \\ 
 \beta_2 = m_{13 }^{2 }= m_{32 }^{1 } =\zeta_1, \quad 
& \alpha_3 = m_{21 }^{4 }= m_{13 }^{2 } =\beta_2,  \\
  \alpha_4= m_{12 }^{4 }= m_{23 }^{1 } =\epsilon_1, \quad 
& \epsilon_1 = m_{23 }^{1 }= m_{31 }^{ 2} =\gamma_2. 
\end{align*}
\end{lemma}
{\bf Proof.}
By \eqref{rs} and the latter half of \eqref{eq017}, 
we have
$r_2 = s_1$. 
From the former half of \eqref{eq017} 
and the coefficients
$(\delta_{c_i} \circ \delta_{c_j}) \circ \delta_{c_k}|_{c_0} 
= \delta_{c_i} \circ (\delta_{c_j} \circ \delta_{c_k})|_{c_0}$,
we have $m_{ij}^{\bar{k}} = m_{jk}^{\bar{i}}$ 
for all $1\le i ,j ,k \le 4$.
The desired results follow.
\hspace{\fill}$\Box$

\bigskip

From Lemma \ref{lem4.5}, 
we can reduce the number of independent parameters as follows:
\begin{align*}
& p = p_0 =q_0 = r_0, \\
& r = r_1 = s_2 = \beta_3 = \zeta_4, \\
& s = r_2 = s_1 =\gamma_4 = \epsilon_3, \\
& t = p_2 = p_3 = \alpha_1 = \beta_1 = \gamma_1,\\
& q = q_1 = q_3 =\alpha_2 = \epsilon_2 = \zeta_2, \\
& u = \beta_4 = \eta_1 = \eta_2 = \eta_4 = \gamma_3 = \epsilon_4 = \zeta_3, \\
& v = \alpha_3 = \beta_2 = \zeta_1, \\
& w = \alpha_4 = \gamma_2 = \epsilon_1,\\
& x= r_3 = s_3 = \eta_3. 
\end{align*}
By axiom (b),
the coefficients of $\delta_{c_1} \circ \delta_{c_3}$ 
and $\delta_{c_2} \circ \delta_{c_4}$ imply that
$\beta_1 + \beta_2 + \beta_3 + \beta_4 =1$
and 
$\zeta_1 + \zeta_2 + \zeta_3 + \zeta_4 =1$,
and hence,
$$
t + v + r + u = 1~~{\rm and}~~v + q + u + r = 1. 
$$
Combining these yields
\begin{equation}
\label{eq018}
t = q.
\end{equation} 
Moreover, the coefficients of
$\delta_{c_1} \circ \delta_{c_1}$ and 
$\delta_{c_2} \circ \delta_{c_2}$
yield the results that
$p_0 + p_1 + p_2 + 2p_3 =1$
and
$q_0 + q_1 + q_2 + 2q_3 =1$;
hence,
$$
p + p_1 + 3t =1~~{\rm and}~~p + q_2 + 3q =1. 
$$
Substituting \eqref{eq018} 
into the above yields
$p_1 = q_2$.
Hence, we can define a parameter $y$ as
$$
y = p_1 = q_2. 
$$

\begin{lemma}
\label{lem4.6}
Let $p, q, r, s, u, v, w, x$ and $y$  be as above.
Then, $p, q, u, v, w, x$ and $y$ are written 
in terms of $r$ and $s$ as follows:
\begin{align*}
& p = 2-4r-4s,~~~q = \frac{1}{2}(r + s),~~~
	u = \frac{1}{2}(1- r - s),\\
& v =  \frac{1}{2}-r,~~w = \frac{1}{2}-s, ~~
	x = \frac{1}{2}(-1+3r+3s),
 	~~y=-1 + \frac{5}{2}(r + s).
\end{align*}
\end{lemma}
{\bf Proof.}
Using the parameters, 
we can write the structure equations of $\mathcal{K}$ as
\begin{align*}
&\delta_{c_1} \circ \delta_{c_1} 
= p \delta_{c_0} + y \delta_{c_1} + q \delta_{c_2} + q \delta_{c_3} + q \delta_{c_4},
\end{align*}
\begin{align*}
&\delta_{c_2} \circ \delta_{c_2} 
= p \delta_{c_0} + q \delta_{c_1} + y \delta_{c_2} + q \delta_{c_3} + q \delta_{c_4},
\end{align*}
\begin{align*}
&\delta_{c_3} \circ \delta_{c_4} 
= p \delta_{c_0} + r \delta_{c_1} + s \delta_{c_2} + x \delta_{c_3} + x \delta_{c_4},\\
&\delta_{c_4} \circ \delta_{c_3} 
= p \delta_{c_0} + s \delta_{c_1} + r \delta_{c_2} + x \delta_{c_3} + x \delta_{c_4},
\end{align*}
\begin{align*}
&\delta_{c_1} \circ \delta_{c_2} 
= q \delta_{c_1} + q \delta_{c_2} + v \delta_{c_3} + w \delta_{c_4},\\
&\delta_{c_2} \circ \delta_{c_1} 
= q \delta_{c_1} + q \delta_{c_2} + w \delta_{c_3} + v \delta_{c_4},
\end{align*}
\begin{align*}
&\delta_{c_1} \circ \delta_{c_3} 
= q \delta_{c_1} + v \delta_{c_2} + r \delta_{c_3} + u \delta_{c_4},\\
&\delta_{c_4} \circ \delta_{c_1} 
= q \delta_{c_1} + v \delta_{c_2} + u \delta_{c_3} + r \delta_{c_4},\\
&\delta_{c_3} \circ \delta_{c_1} 
= q \delta_{c_1} + w \delta_{c_2} + s \delta_{c_3} + u \delta_{c_4},\\
&\delta_{c_1} \circ \delta_{c_4} 
= q \delta_{c_1} + w \delta_{c_2} + u \delta_{c_3} + s \delta_{c_4},
\end{align*}
\begin{align*}
&\delta_{c_2} \circ \delta_{c_3} 
= w \delta_{c_1} + q \delta_{c_2} + s \delta_{c_3} + u \delta_{c_4},\\
&\delta_{c_4} \circ \delta_{c_2} 
= w \delta_{c_1} + q \delta_{c_2} + u \delta_{c_3} + s \delta_{c_4},\\
&\delta_{c_3} \circ \delta_{c_2} 
= v \delta_{c_1} + q \delta_{c_2} + r \delta_{c_3} + u \delta_{c_4},\\
&\delta_{c_2} \circ \delta_{c_4} 
= v \delta_{c_1} + q \delta_{c_2} + u \delta_{c_3} + r \delta_{c_4},
\end{align*}
\begin{align*}
&\delta_{c_3} \circ \delta_{c_3} 
= u \delta_{c_1} + u  \delta_{c_2} + x \delta_{c_3} + u \delta_{c_4},\\
&\delta_{c_4} \circ \delta_{c_4} 
= u \delta_{c_1} + u  \delta_{c_2} + u \delta_{c_3} + x \delta_{c_4}. 
\end{align*}
Hence, by axiom (b), 
\begin{align*}
p + y + 3q = 1,~~~&p + r + s + 2x = 1,~~~2q + v + w = 1,\\
q + v + r + u =1,~~~&q + w + u + s = 1,~~~3u + x = 1. 
\end{align*}
With these equations,
the parameters $u$, $x$, $p$, $v$, $w$, and $y$ can be expressed 
in terms of $r$, $s$, and $q$ as follows:
\begin{align*}
u = \frac{1}{2}(1- r - s),~~~&x = \frac{1}{2}(-1+3r+3s),~~~p = 2-4r-4s,\\
 v = \frac{1}{2}(1-r+s)-q,~~~&w = \frac{1}{2}(1+ r -s) -q,~~~y=-1+4r+4s-3q.
\end{align*}
It suffices to express $q$ in terms of $r$ and $s$.
Calculating the coefficients
$(\delta_{c_3} \circ \delta_{c_4}) \circ \delta_{c_3}|_{c_1} 
=\delta_{c_3} \circ (\delta_{c_4} \circ \delta_{c_3})|_{c_1}$
yields
\begin{align*}
rq + sw + xs = sq + rv + xr.
\end{align*}
Because $rv = rw-r^2+rs$, 
$$
q(r-s) = (r-s)(w-r+x). 
$$
By assumption, 
$r \not = s$, and hence,
\begin{align*}
q = w-r +x = \frac{1}{2}(r + s).
\end{align*}
From this, we obtain the expressions for $v$, $w$, and $y$.
This completes the proof.
\hspace{\fill}$\Box$

\bigskip

By axiom (b), 
each coefficient of the structure equations of 
$\mathcal{K}$ is non-negative. Moreover, by axiom (c), $p$ must be positive.
Hence, we get the ranges of the parameters as follows:
\begin{align*}
 r+s \leq 1 \quad {\rm by} \quad 0 \leq u, \qquad
&\frac{1}{3} \leq r+s \quad {\rm by}\quad 0 \leq x,\\
 r+s < \frac{1}{2}\quad {\rm by} \quad 0 < p , \qquad
& r \leq \frac{1}{2}\quad {\rm by}\quad 0 \leq v,\\
 s \leq \frac{1}{2}\quad {\rm by} \quad 0 \leq w , \qquad
&\frac{2}{5} \leq r+s \quad {\rm by}\quad 0 \leq y ,\\
0 \leq r+s \quad {\rm by}\quad 0 \leq q . \qquad &
\end{align*}
Taking the above into consideration, 
we know that if $r$ and $s$ satisfy
\begin{equation}
\label{eq021}
0 \leq r < \frac{1}{2},~~~0 \leq s < \frac{1}{2}~~{\rm and}~~
\frac{2}{5} \leq r + s < \frac{1}{2},
\end{equation}
then $\mathcal{K}$ satisfies conditions (b) and (c) of the axiom.

Finally, we check the associativity of $\mathcal{K}$. 
Since the structure of $\mathcal{K}$ is determined by $r$ and $s$,
we write $\mathcal{K} = \mathcal{K}_{(r,s)}$.
Interchanging the parameters $r$ and $s$,
we set $\mathcal{K}^\prime = \mathcal{K}_{(s,r)}$.
To avoid confusion,
we mark the parameters of $\mathcal{K}^\prime$ with primes.
Then, we observe from Lemma \ref{lem4.6}
that 
\begin{eqnarray*}
&&p=p', \quad q = q', \quad r =s', \quad s = r',  \quad u=u', \\
&& v = w', \quad w = v',
\quad x = x', \quad y=y'. 
\end{eqnarray*}
Define a linear map $\pi$ from 
$M^b(K)$ 
to $M^b(K^\prime)$ 
as
\[
\pi(\delta_{c_0}) = \delta_{c_0} ,\quad 
\pi(\delta_{c_1}) = \delta_{c_2} ,\quad 
\pi(\delta_{c_2}) = \delta_{c_1} ,\quad 
\pi(\delta_{c_3}) = \delta_{c_3} ,\quad 
\pi(\delta_{c_4}) = \delta_{c_4}.
\]
\begin{lemma}
\label{lem4.7}
$\pi$ preserves the convolution $\circ$.
\end{lemma}
{\bf Proof.}
The equation
$
\delta_{c_1} \circ \delta_{c_1} 
= p \delta_{c_0} + y \delta_{c_1} + q \delta_{c_2} + q \delta_{c_3} + q \delta_{c_4}
$
implies that
\[
\pi(\delta_{c_1} \circ \delta_{c_1})
= p \delta_{c_0} + q \delta_{c_1} + y \delta_{c_2} + q \delta_{c_3} + q \delta_{c_4}.
\]
On the other hand,
\[
\pi(\delta_{c_1}) \circ \pi(\delta_{c_1})
=\delta_{c_2} \circ \delta_{c_2} 
=p' \delta_{c_0} + q' \delta_{c_1} + y' \delta_{c_2} + q' \delta_{c_3} + q' \delta_{c_4}.
\]
Therefore, 
$\pi(\delta_{c_1} \circ \delta_{c_1}) = \pi(\delta_{c_1}) \circ \pi(\delta_{c_1})$.
Similarly, 
\[
\pi( \delta_{c_1} \circ \delta_{c_3})
= v \delta_{c_1} + q \delta_{c_2} + r \delta_{c_3} + u \delta_{c_4}
\]
and
\[
\pi(\delta_{c_1} )\circ \pi( \delta_{c_3})
=\delta_{c_2} \circ \delta_{c_3}
= w' \delta_{c_1} + q' \delta_{c_2} + s' \delta_{c_3} + u' \delta_{c_4}.
\]
Hence,  
$\pi(\delta_{c_1} \circ \delta_{c_3}) = \pi(\delta_{c_1}) \circ \pi(\delta_{c_3})$.
The others are left to the reader.
\hspace{\fill}$\Box$

\bigskip
 
Lemma \ref{lem4.7} says that if the associativity
\[
(\delta_{c_i} \circ \delta_{c_j}) \circ \delta_{c_k} 
= \delta_{c_i} \circ (\delta_{c_j} \circ \delta_{c_k})
\]
holds for all $r$ and $s$, then the associativity 
\[
(\pi(\delta_{c_i}) \circ \pi(\delta_{c_j})) \circ \pi(\delta_{c_k} )
= \pi(\delta_{c_i}) \circ (\pi(\delta_{c_j}) \circ \pi(\delta_{c_k}))
\]
also holds.

Let $\tilde \pi$ be a linear map
 from $M^b(K)$ to $M^b(K^\prime)$ defined as
\[
\tilde\pi(\delta_{c_0}) = \delta_{c_0} ,\quad 
\tilde\pi(\delta_{c_1}) = \delta_{c_1} ,\quad 
\tilde\pi(\delta_{c_2}) = \delta_{c_2} ,\quad 
\tilde\pi(\delta_{c_3}) = \delta_{c_4} ,\quad 
\tilde\pi(\delta_{c_4}) = \delta_{c_3}.
\]
In a manner similar to the proof of Lemma \ref{lem4.7},
we can prove the following.
\begin{lemma}
\label{lem4.8}
$\tilde\pi$
preserves the involution $\circ$.
\end{lemma}
\begin{theorem}
\label{thm4.9}
Let $r$ and $s$ be as above.
Suppose that \eqref{eq021}
and
\begin{equation}
\label{eq020}
3r^2 + 10rs + 3s^2 - 8r -8s +3 = 0.
\end{equation}
Then, $\mathcal{K}$ is a non-commutative hypergroup.
\end{theorem}
\noindent
{\bf Proof.}
From the above argument, 
\eqref{eq021} implies that $\mathcal{K} = \mathcal{K}_{(r,s)}$ satisfies 
conditions (b) and (c) of the axiom.
It suffices to check the associativities.
It is clear that if  
$(\delta_{c_i} \circ \delta_{c_j}) \circ \delta_{c_k} 
= \delta_{c_i} \circ (\delta_{c_j} \circ \delta_{c_k})$, 
then 
$(\delta_{c_{\bar{k}}} \circ \delta_{c_{\bar{j}}}) \circ \delta_{c_{\bar{i}}} 
= \delta_{c_{\bar{k}}} \circ (\delta_{c_{\bar{j}}} \circ \delta_{c_{\bar{i}}})$ 
for all $1 \le i,j,k \le 4$.
Combining this with Lemmas \ref{lem4.7} and \ref{lem4.8}, 
we only need to check the following:  
\begin{equation}
\label{eq022}
\begin{cases}
(\delta_{c_1} \circ \delta_{c_1}) \circ \delta_{c_1} 
= \delta_{c_1} \circ (\delta_{c_1} \circ \delta_{c_1}),~~~
(\delta_{c_1} \circ \delta_{c_1}) \circ \delta_{c_2} 
= \delta_{c_1} \circ (\delta_{c_1} \circ \delta_{c_2}),\\
(\delta_{c_1} \circ \delta_{c_1}) \circ \delta_{c_3} 
= \delta_{c_1} \circ (\delta_{c_1} \circ \delta_{c_3}),~~~
(\delta_{c_1} \circ \delta_{c_2}) \circ \delta_{c_1} 
= \delta_{c_1} \circ (\delta_{c_2} \circ \delta_{c_1}),\\
(\delta_{c_1} \circ \delta_{c_2}) \circ \delta_{c_3} 
= \delta_{c_1} \circ (\delta_{c_2} \circ \delta_{c_3}),~~~
(\delta_{c_1} \circ \delta_{c_3}) \circ \delta_{c_1} 
= \delta_{c_1} \circ (\delta_{c_3} \circ \delta_{c_1}),\\
(\delta_{c_1} \circ \delta_{c_3}) \circ \delta_{c_2} 
= \delta_{c_1} \circ (\delta_{c_3} \circ \delta_{c_2}),~~~
(\delta_{c_1} \circ \delta_{c_3}) \circ \delta_{c_3} 
= \delta_{c_1} \circ (\delta_{c_3} \circ \delta_{c_3}),\\
(\delta_{c_1} \circ \delta_{c_3}) \circ \delta_{c_4} 
= \delta_{c_1} \circ (\delta_{c_3} \circ \delta_{c_4}),~~~
(\delta_{c_3} \circ \delta_{c_1}) \circ \delta_{c_3} 
= \delta_{c_3} \circ (\delta_{c_1} \circ \delta_{c_3}),\\
(\delta_{c_3} \circ \delta_{c_1}) \circ \delta_{c_4} 
= \delta_{c_3} \circ (\delta_{c_1} \circ \delta_{c_4}),~~~
(\delta_{c_3} \circ \delta_{c_3}) \circ \delta_{c_3} 
= \delta_{c_3} \circ (\delta_{c_3} \circ \delta_{c_3}),\\
(\delta_{c_3} \circ \delta_{c_3}) \circ \delta_{c_4} 
= \delta_{c_3} \circ (\delta_{c_3} \circ \delta_{c_4}),~~~
(\delta_{c_3} \circ \delta_{c_4}) \circ \delta_{c_3} 
= \delta_{c_3} \circ (\delta_{c_4} \circ \delta_{c_3}).
\end{cases}
\end{equation}
We provide the proof only for the first two equations in \eqref{eq022};
the other equations are proved in \ref{appe1}.
Direct calculation yields
\begin{align*}
(\delta_{c_1} \circ \delta_{c_1}) \circ \delta_{c_1}  
&= yp \delta_{c_0} + (p + y^2 + 3q^2)\delta_{c_1} + (yq + q^2 + qw + qv)\delta_{c_2} \\
&+(yq + qw + qs + qu)\delta_{c_3} + (yq + qv + qu + qr)\delta_{c_4}, \\[3mm]
\delta_{c_1} \circ (\delta_{c_1} \circ \delta_{c_1})
&= yp \delta_{c_0} + (p + y^2 + 3q^2)\delta_{c_1} + (yq + q^2 + qv + qw)\delta_{c_2} \\
&+(yq + qv + qr + qu)\delta_{c_3} + (yq + qw + qu + qs)\delta_{c_4}.
\end{align*}
Hence, by Lemma \ref{lem4.6},
we have $(\delta_{c_1} \circ \delta_{c_1}) \circ \delta_{c_1} 
= \delta_{c_1} \circ (\delta_{c_1} \circ \delta_{c_1})$.

To show that $(\delta_{c_1} \circ \delta_{c_1}) \circ \delta_{c_2} 
= \delta_{c_1} \circ (\delta_{c_1} \circ \delta_{c_2})$,
it suffices to prove that
\begin{equation}
\label{eq019}
\begin{cases}
p + 2yq = v^2 + w^2, \\
yv + qr + qu = qv + vr + wu,\\
yw + qu + qs = qw + vu + ws,
\end{cases}
\end{equation}
because
\begin{align*}
(\delta_{c_1} \circ \delta_{c_1}) \circ \delta_{c_2}  
&= qp \delta_{c_0} + (yq + q^2 + qv + qw)\delta_{c_1} + (p + 2yq + 2q^2)\delta_{c_2} \\
&+(yv + q^2 + qr + qu)\delta_{c_3} + (yw + q^2 + qu + qs)\delta_{c_4}, \\[3mm]
\delta_{c_1} \circ (\delta_{c_1} \circ \delta_{c_2})
&= qp \delta_{c_0} + (qy + q^2 + vq + wq )\delta_{c_1} + (2q^2 + v^2 + w^2)\delta_{c_2} \\
&+(q^2 + qv + vr + wu)\delta_{c_3} + (q^2 + qw + vu + ws)\delta_{c_4}.
\end{align*}
We observe from Lemma \ref{lem4.6} that \eqref{eq020} 
yields \eqref{eq019}.
\hspace{\fill}$\Box$

\bigskip

Let $\Lambda$ be the set of all pairs $(r,s)$
that satisfy \eqref{eq021} and \eqref{eq020}.
$\Lambda$ is not empty (see the figure below).
Let $(r,s) \in \Lambda$. 
By Theorem \ref{thm4.9}, 
$\mathcal{K} = \mathcal{K}_{(r,s)}$ is a non-commutative hypergroup.
This proves Theorem \ref{thm}.

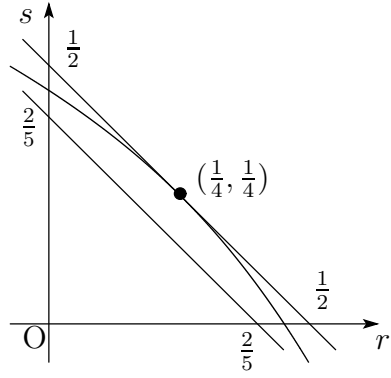
\begin{center}
\input{non_commutative_zu1.tex}
\end{center}


\appendix
\def\thesection{Appendix \Alph{section}}
\section{}\label{appe1}

In this appendix, we show the details of the calculation of associativity. 
Since $(\delta_{c_1} \circ \delta_{c_1}) \circ \delta_{c_1} 
= \delta_{c_1} \circ (\delta_{c_1} \circ \delta_{c_1})$
and  $(\delta_{c_1} \circ \delta_{c_1}) \circ \delta_{c_2} 
= \delta_{c_1} \circ (\delta_{c_1} \circ \delta_{c_2})$
are proved in the proof of Theorem \ref{thm4.9}, we only need to check
 the equations in \eqref{eq022}
 other than $(\delta_{c_1} \circ \delta_{c_1}) \circ \delta_{c_1} 
= \delta_{c_1} \circ (\delta_{c_1} \circ \delta_{c_1})$
and  $(\delta_{c_1} \circ \delta_{c_1}) \circ \delta_{c_2} 
= \delta_{c_1} \circ (\delta_{c_1} \circ \delta_{c_2})$.

By the calculations
\begin{align*}
(\delta_{c_1} \circ \delta_{c_1}) \circ \delta_{c_3}  
&= qp \delta_{c_0} + (yq + qw + qu + qs)\delta_{c_1} + (yv + q^2 + qu + qr)\delta_{c_2} \\
&+(p  + yr + qs + 2qx)\delta_{c_3} + (yu + 2qu + qx)\delta_{c_4}, \\[3mm]
\delta_{c_1} \circ (\delta_{c_1} \circ \delta_{c_3})
&= qp \delta_{c_0} + (qy + vq + rq + uq )\delta_{c_1} + (q^2 + vq + rv + uw)\delta_{c_2} \\
&+(q^2 + v^2 + r^2 + u^2)\delta_{c_3} + (q^2 + vw + ru + us)\delta_{c_4}, 
\end{align*}
we need
\begin{align*}
qw + qs &= vq + rq,\\
yv +  qu + qr &= vq + rv + uw,\\
p + yr + qs + 2qx &= q^2 + v^2 + r^2 + u^2,\\
yu + 2qu +  qx &= q^2 + vw + ru + us 
\end{align*}
for $(\delta_{c_1} \circ \delta_{c_1}) \circ \delta_{c_3} 
= \delta_{c_1} \circ (\delta_{c_1} \circ \delta_{c_3})$.
The first equation always holds. The equation
\begin{equation}\label{conditionrs}
3r^2 + 10rs + 3s^2 - 8r -8s +3 = 0
\end{equation}
is required to satisfy the others.


By the calculations
\begin{align*}
(\delta_{c_1} \circ \delta_{c_2}) \circ \delta_{c_1}  
&= qp \delta_{c_0} + (qy + q^2 + vq + wq)\delta_{c_1} + (2q^2 + 2vw)\delta_{c_2} \\
&+(q^2  + qw + vs + wu)\delta_{c_3} + (q^2 + qv + vu + wr)\delta_{c_4}, \\[3mm]
\delta_{c_1} \circ (\delta_{c_2} \circ \delta_{c_1})
&= qp \delta_{c_0} + (qy + q^2 + wq + vq )\delta_{c_1} + (2q^2 + 2vw)\delta_{c_2} \\
&+(q^2 + qv + wr + vu)\delta_{c_3} + (q^2 + qw + wu + vs)\delta_{c_4},
\end{align*}
we need 
\begin{align*}
qw + vs + wu &= qv + wr + vu 
\end{align*}
for $(\delta_{c_1} \circ \delta_{c_2}) \circ \delta_{c_1} 
= \delta_{c_1} \circ (\delta_{c_2} \circ \delta_{c_1})$. 
This equation always holds. 


By the calculations
\begin{align*}
(\delta_{c_1} \circ \delta_{c_2}) \circ \delta_{c_3}  
&= wp \delta_{c_0} + (q^2 + qw + vu + ws)\delta_{c_1} + (qv + q^2 + vu + wr)\delta_{c_2} \\
&+(qr + qs + vx + wx)\delta_{c_3} + (2qu + vu + wx)\delta_{c_4}, \\[3mm]
\delta_{c_1} \circ (\delta_{c_2} \circ \delta_{c_3})
&= wp \delta_{c_0} + (wy + q^2 + sq + uq )\delta_{c_1} + (wq + q^2 + sv + uw)\delta_{c_2} \\
&+(wq + qv + sr + u^2)\delta_{c_3} + (2wq  + 2su)\delta_{c_4},
\end{align*}
we need  
\begin{align*}
qw + vu + ws &= wy + sq + uq,\\
qv + vu + wr &= wq + sv + uw,\\
qr + qs + vx + wx &= wq + qv + sr + u^2,\\
2qu + vu + wx &= 2wq + 2su 
\end{align*}
for $(\delta_{c_1} \circ \delta_{c_2}) \circ \delta_{c_3} 
= \delta_{c_1} \circ (\delta_{c_2} \circ \delta_{c_3})$. 
The second and fourth equations always hold.  The equation \eqref{conditionrs}
is required to satisfy the others.

By the calculations
\begin{align*}
(\delta_{c_1} \circ \delta_{c_3}) \circ \delta_{c_1}  
&= qp \delta_{c_0} + (qy + vq + rq + uq)\delta_{c_1} + (q^2 + vq + rw + uv)\delta_{c_2} \\
&+(q^2 + vw + rs + u^2)\delta_{c_3} + (q^2 + v^2 + 2ru)\delta_{c_4}, \\[3mm]
\delta_{c_1} \circ (\delta_{c_3} \circ \delta_{c_1})
&= qp \delta_{c_0} + (qy + wq + sq + uq )\delta_{c_1} + (q^2 + wq + sv + uw)\delta_{c_2} \\
&+(q^2 + wv + sr + u^2)\delta_{c_3} + (q^2 + w^2 + 2su)\delta_{c_4},
\end{align*}
we need 
\begin{align*}
vq + rq  &= wq + sq,\\
vq + rw + uv &= wq + sv + uw,\\
v^2 + 2ru &= w^2 + 2su 
\end{align*}
for $(\delta_{c_1} \circ \delta_{c_3}) \circ \delta_{c_1} 
= \delta_{c_1} \circ (\delta_{c_3} \circ \delta_{c_1})$.
These equations always hold.

By the calculations
\begin{align*}
(\delta_{c_1} \circ \delta_{c_3}) \circ \delta_{c_2}  
&= vp \delta_{c_0} + (q^2 + vq + rv + uw)\delta_{c_1} + (q^2 + vy + rq + uq)\delta_{c_2} \\
&+(2qv + r^2 + u^2)\delta_{c_3} + (qw + vq + ru + us)\delta_{c_4}, \\[3mm]
\delta_{c_1} \circ (\delta_{c_3} \circ \delta_{c_2})
&= vp \delta_{c_0} + (vy + q^2 + rq + uq )\delta_{c_1} + (vq + q^2 + rv + uw)\delta_{c_2} \\
&+(2vq  + r^2 + u^2)\delta_{c_3} + (vq + qw + ru + us)\delta_{c_4},
\end{align*}
we need 
\begin{align*}
vq + rv + uw  &= vy + rq + uq,
\end{align*}
for $(\delta_{c_1} \circ \delta_{c_3}) \circ \delta_{c_2} 
= \delta_{c_1} \circ (\delta_{c_3} \circ \delta_{c_2})$.
The equation \eqref{conditionrs}
is required to satisfy the above equation.


By the calculations
\begin{align*}
(\delta_{c_1} \circ \delta_{c_3}) \circ \delta_{c_3}  
&= up \delta_{c_0} + (q^2 + vw + ru + us)\delta_{c_1} + (2qv  + 2ru)\delta_{c_2} \\
&+(qr + vs + rx + ux)\delta_{c_3} + (qu + vu + ru + ux)\delta_{c_4}, \\[3mm]
\delta_{c_1} \circ (\delta_{c_3} \circ \delta_{c_3})
&= up \delta_{c_0} + (uy + 2uq + xq)\delta_{c_1} + (2uq  + xv + uw)\delta_{c_2} \\
&+(uq + uv + xr + u^2)\delta_{c_3} + (uq + uw + xu + us)\delta_{c_4},
\end{align*}
we need 
\begin{align*}
q^2 + vw + ru + us &= uy + 2uq + xq,\\
2qv  + 2ru &= 2uq + xv + uw,\\
qr + vs  + ux &= uq + uv  + u^2,\\
vu + ru  &=  uw  + us
\end{align*}
for $(\delta_{c_1} \circ \delta_{c_3}) \circ \delta_{c_3} 
= \delta_{c_1} \circ (\delta_{c_3} \circ \delta_{c_3})$. 
The second and fourth equations always hold.
The equation \eqref{conditionrs}
is required to satisfy the others.

By the calculations
\begin{align*}
(\delta_{c_1} \circ \delta_{c_3}) \circ \delta_{c_4}  
&= rp \delta_{c_0} + (q^2 + v^2 + r^2 + u^2)\delta_{c_1} + (qw + vq + rs + u^2)\delta_{c_2} \\
&+(qu + vu + rx + u^2)\delta_{c_3} + (qs + vr + rx + ux)\delta_{c_4}, \\[3mm]
\delta_{c_1} \circ (\delta_{c_3} \circ \delta_{c_4})
&= rp \delta_{c_0} + (p + ry + sq + 2xq)\delta_{c_1} + (rq + sq + xv + xw)\delta_{c_2} \\
&+(rq + sv + xr + xu)\delta_{c_3} + (rq + sw + xu + xs)\delta_{c_4},
\end{align*}
we need  
\begin{align*}
q^2 + v^2 + r^2 + u^2 &= p + ry + sq + 2xq,\\
qw + vq + rs + u^2 &= rq + sq + xv + xw,\\
qu + vu  + u^2 &= rq + sv + xu,\\
qs + vr + rx  &=  rq + sw + xs
\end{align*}
for $(\delta_{c_1} \circ \delta_{c_3}) \circ \delta_{c_4} 
= \delta_{c_1} \circ (\delta_{c_3} \circ \delta_{c_4})$.
The fourth equation always holds.
The equation \eqref{conditionrs}
is required to satisfy the others.


By the calculations
\begin{align*}
(\delta_{c_3} \circ \delta_{c_1}) \circ \delta_{c_3}  
&= up \delta_{c_0} + (q^2 + w^2 + 2su)\delta_{c_1} + (qv + wq + su + ur)\delta_{c_2} \\
&+(qr + ws + sx + ux)\delta_{c_3} + (qu + wu + su + ux)\delta_{c_4}, \\[3mm]
\delta_{c_3} \circ (\delta_{c_1} \circ \delta_{c_3})
&= up \delta_{c_0} + (q^2 + v^2 + 2ru)\delta_{c_1} + (qw + vq + ru + us)\delta_{c_2} \\
&+(qs + vr + rx + ux)\delta_{c_3} + (qu + vu + ru + ux)\delta_{c_4},
\end{align*}
we need 
\begin{align*}
w^2 + 2su &= v^2 + 2ru,\\
qr + ws  + sx &= qs + vr + rx,\\
wu + su  &=  vu + ru 
\end{align*}
for $(\delta_{c_3} \circ \delta_{c_1}) \circ \delta_{c_3} 
= \delta_{c_3} \circ (\delta_{c_1} \circ \delta_{c_3})$. 
These equations always hold.

By the calculations
\begin{align*}
(\delta_{c_3} \circ \delta_{c_1}) \circ \delta_{c_4}  
&= sp \delta_{c_0} + (q^2 + wv + sr + u^2)\delta_{c_1} + (2qw  + s^2 + u^2)\delta_{c_2} \\
&+(qu + wu + sx + u^2)\delta_{c_3} + (qs + wr + sx + ux)\delta_{c_4}, \\[3mm]
\delta_{c_3} \circ (\delta_{c_1} \circ \delta_{c_4})
&= sp \delta_{c_0} + (q^2 + wv + u^2 + sr)\delta_{c_1} + (2qw  + u^2 + s^2)\delta_{c_2} \\
&+(qs + wr + ux + sx)\delta_{c_3} + (qu + wu + u^2 + sx)\delta_{c_4},
\end{align*}
we need 
\begin{align*}
qu + wu + u^2 &= qs + wr + ux
\end{align*}
for $(\delta_{c_3} \circ \delta_{c_1}) \circ \delta_{c_4} 
= \delta_{c_3} \circ (\delta_{c_1} \circ \delta_{c_4})$.
The equation \eqref{conditionrs}
is required to satisfy the above equation.

By the calculations
\begin{align*}
(\delta_{c_3} \circ \delta_{c_3}) \circ \delta_{c_3}  
&= up \delta_{c_0} + (uq + uw + xu + us)\delta_{c_1} + (uv + uq + xu + ur)\delta_{c_2} \\
&+(ur + us + x^2 + ux)\delta_{c_3} + (2u^2 + 2ux)\delta_{c_4}, \\[3mm]
\delta_{c_3} \circ (\delta_{c_3} \circ \delta_{c_3})
&= up \delta_{c_0} + (uq + uv + xu + ur)\delta_{c_1} + (uw + uq + xu + us)\delta_{c_2} \\
&+(us + ur + x^2 + ux)\delta_{c_3} + (2u^2 + 2xu)\delta_{c_4},
\end{align*}
we need 
\begin{align*}
uw + us  &= uv + ur
\end{align*}
for $(\delta_{c_3} \circ \delta_{c_3}) \circ \delta_{c_3} 
= \delta_{c_3} \circ (\delta_{c_3} \circ \delta_{c_3})$.
This equation always holds.


By the calculations
\begin{align*}
(\delta_{c_3} \circ \delta_{c_3}) \circ \delta_{c_4}  
&= xp \delta_{c_0} + (uq + uv + xr + u^2)\delta_{c_1} + (uw + uq + xs + u^2)\delta_{c_2} \\
&+(3u^2 + x^2)\delta_{c_3} + (us + ur + x^2 + ux)\delta_{c_4}, \\[3mm]
\delta_{c_3} \circ (\delta_{c_3} \circ \delta_{c_4})
&= xp \delta_{c_0} + (rq + sv + xu + xr)\delta_{c_1} + (rw + sq + xu + xs)\delta_{c_2} \\
&+(p + 2rs  + 2x^2)\delta_{c_3} + (ru + su + xu + x^2)\delta_{c_4},
\end{align*}
we need  
\begin{align*}
uq + uv + u^2  &= rq + sv + xu,\\
uw + uq + u^2 &= rw + sq + xu,\\
3u^2  &= p + 2rs + x^2
\end{align*}
for $(\delta_{c_3} \circ \delta_{c_3}) \circ \delta_{c_4} 
= \delta_{c_3} \circ (\delta_{c_3} \circ \delta_{c_4})$.
The equation \eqref{conditionrs}
is required to satisfy the above equations.

By the calculations
\begin{align*}
(\delta_{c_3} \circ \delta_{c_4}) \circ \delta_{c_3}  
&= xp \delta_{c_0} + (rq + sw + xu + xs)\delta_{c_1} + (rv + sq + xu + xr)\delta_{c_2} \\
&+(p + r^2 + s^2 + 2x^2)\delta_{c_3} + (ru + su + xu + x^2)\delta_{c_4}, \\[3mm]
\delta_{c_3} \circ (\delta_{c_4} \circ \delta_{c_3})
&= xp \delta_{c_0} + (sq + rv + xu + xr)\delta_{c_1} + (sw + rq + xu + xs)\delta_{c_2} \\
&+(p + s^2 + r^2 + 2x^2)\delta_{c_3} + (su + ru + xu + x^2)\delta_{c_4}, 
\end{align*}
we need  
\begin{align*}
rq + sw + xs &= sq + rv + xr
\end{align*}
for $(\delta_{c_3} \circ \delta_{c_4}) \circ \delta_{c_3} 
= \delta_{c_3} \circ (\delta_{c_4} \circ \delta_{c_3})$.
This equation always holds.

\bigskip
\noindent
{\bf Acknowledgement.} This work was supported by 
JSPS KAKENHI Grant Numbers 25800061, 26350231, 26800055.

\bigskip

\textbf{Addresses }

\medskip

Yasumichi Matsuzawa : Department of Mathematics, Faculty of Education, Shinshu University,
6-Ro, Nishi-nagano, Nagano, 380-8544, Japan

\medskip

e-mail : myasu@shinshu-u.ac.jp

\bigskip

\medskip

Hiromichi Ohno : Department of Mathematics, Faculty of Engineering, Shinshu University,
4-17-1 Wakasato, Nagano, 380-8553, Japan

\medskip

e-mail : h\_ohno@shinshu-u.ac.jp

\bigskip

\medskip

Akito Suzuki : Department of Mathematics, Faculty of Engineering, Shinshu University,
4-17-1 Wakasato, Nagano, 380-8553, Japan

\medskip

e-mail : akito@shinshu-u.ac.jp

\bigskip

Tatsuya Tsurii : Graduate School of Science of the Osaka Prefecture University 
1-1 Gakuen-cho, Nakaku, Sakai,
Osaka, 599-8531, Japan

\medskip

e-mail : dw301003@edu.osakafu-u.ac.jp

\bigskip

Satoe Yamanaka : 
Faculty of science,
Nara Women's University,
Kitauoyahigashimachi, Nara, 630-8506, Japan

\medskip

e-mail : s.yamanaka@gmail.com

\end{document}